# A zero-one law for dynamical properties

Eli Glasner & Jonathan L. King

*Tel Aviv University, Tel Aviv, Israel*   `glasner@math.tau.ac.il`
*University of Florida, Gainesville 32611-2082, USA*   `squash@math.ufl.edu`

ABSTRACT. For any countable group $\Gamma$ satisfying the "weak Rohlin property", and for each dynamical property, the set of $\Gamma$-actions with that property is either residual or meager. The class of groups with the weak Rohlin property includes each lattice $\mathbb{Z}^{\times d}$; indeed, all countable discrete amenable groups.

For $\Gamma$ an arbitrary countable group, let $\mathbb{A}$ be the set of $\Gamma$-actions on the unit circle $Y$. We establish an Equivalence theorem by showing that a dynamical property is Baire/meager/residual in $\mathbb{A}$ if and only if it is Baire/meager/residual in the set of shift-invariant measures on the product space $Y^{\times \Gamma}$.

## §1   INTRODUCTION

Halmos's book *Ergodic Theory* introduced many of us to the study of determining which dynamical properties are ***generic*** (i.e, topologically residual) in the so-called "coarse topology" on transformations. For instance, "weak-mixing" is generic, whereas "mixing" is not, [Hal, pp. 77,78]. The exploration of this notion of genericity became an active research area; see [CP] and [CN] for results and extensive bibliographies.

It has often been the case, when a property has failed to be generic, that further investigation has shown its negation to be generic. This suggests that there is a type of **"zero-one"** law operating for dynamical properties –each is either ***meager*** or ***residual***; that is, either a first-category set (a countable union of nowhere-dense sets) or the complement of such.

In 1993, we found a demonstration of this for certain acting groups $\Gamma$, for dynamical properties which are Baire-measurable. We first present the brief proof in the abstract framework of a group $\Phi$ acting as homeomorphisms of a topological space $\mathbb{A}$, where $\mathbb{A}$ fulfills the Baire Category Theorem in that each residual set is dense. We call such a space a ***BaireCat space***, and say that a subset of a topological space is ***BaireCat*** if it is a BaireCat space in its induced topology.[†]

We then show how this abstract framework applies when $\mathbb{A}$ is the space of measure-preserving $\Gamma$-actions with $\Phi$ playing the role of its group of isomorphisms. In this instance, $\mathbb{A}$ will be a ***Polish space***, that is, homeomorphic to a complete separable metric space.

---

1991 *Mathematics Subject Classification*. 28D15. Secondary: 54E52 60F20 54H05.

*Key words and phrases.* Ergodic theory, Baire Category, Genericity.

J.L. King partially supported by NSF grant DMS-9112595.

[†]A "Baire space" is the traditional name for what we call a BaireCat space. The traditional meanings of "Baire set" and "Baire space" are not related by the induced topology. In view of this, we have revised terminology for this article.

*Miscellany.* In this article, a "field" is a sigma-algebra and "measure" means a probability measure on the Borel-field of a topological space. All sets and mappings are Borel-measurable unless specified otherwise. Given a measure $\mu$, we call $\mu(A)$ the $\mu$-***mass*** of $A$. Let ***mp*** abbreviate "measure preserving".

For a set, measure or transformation $A$, let $A^{\times 3}$ denote the cartesian power $A \times A \times A$. For a set $G$, let $A^{\times G}$ means the cartesian power $\bigotimes_{\gamma \in G} A_\gamma$, where each $A_\gamma = A$.

**What is a Dynamical Property?** Leaving the definition of coarse topology until later, consider the case of a measure-preserving transformation (***mpt***) $T$ of a probability space $\mathbb{Y} = (Y, \mathcal{Y}, \mu)$. Those transformations isomorphic to $T$ are of the form $\varphi T \varphi^{-1}$, where $\varphi$ is a member of the group, $\Phi$, of bi-mpts. [We write "$\varphi(T)$" for $\varphi\, T \varphi^{-1}$, to notationally emphasize that $\Phi$ *acts* on the set of transformations.] Consequently, a dynamical property, such as the set of mixing transformations, is necessarily a union of isomorphism classes, that is, a union of $\Phi$-orbits.

But $\Phi$-invariance is not enough. For an invariant set, $B$, of transformations to be considered a "dynamical property", it is reasonable to require that it satisfy some kind of measurability condition. We can prove the Zero-One Law when $B$ is Borel-measurable or even (co-)analytic, because our proof applies when $B$ is a ***Baire set*** –that is, when $B$ can be written as the symmetric-difference $V \bigtriangleup M$ of an open set with a meager set $M$. Call $V \bigtriangleup M$ a "picture" of $B$; a Baire set may have many pictures. The Baire sets[‡] form a field which includes the Borel sets—indeed, even the analytic sets; see [Kur, pp. 482, 92]. As evidence that dynamical properties might transcend the Borel-field, the Prime Example, below, considers the set of actions without factors.

Before showing its application to transformation-groups, we now state and prove the result in an abstract setting. The abstract lemma, we later found out, appears as "folklore" (private communication) in [Kec, 1995], theorem 8.46. The earliest reference to the result that we have found is in [Oxt2, 1937], where it is stated in the case that $\Phi$ is the group of integers.

Zero-one Lemma. *Suppose $\Phi$ is a group of homeomorphisms of a BaireCat space $\mathbb{A}$. If there is some member $\mathbf{T} \in \mathbb{A}$ whose $\Phi$-orbit is dense[‡] in $\mathbb{A}$, then: Each Baire-measurable $\Phi$-invariant set $B$ is either meager or residual.*

**The Baire necessities.** (Proofs of the following facts can be found in [Oxt, Chap. 4].) A subset of a topological space is ***regular-open*** if it equals the interior of its closure. The fact we need from elementary topology is that a Baire set $B$ has a "regular-picture" $B = U \bigtriangleup M$, where $U$ is regular-open and $M$ is meager.

The following four properties of a topological space $\mathbb{A}$ are easily equivalent: *(a) Each Baire set $B \subset \mathbb{A}$ has at most one (hence exactly one) regular-picture $B = U \bigtriangleup M$. (b/c) The only open/regular-open meager set is empty. (d) $\mathbb{A}$ is BaireCat.* (By the way, the Baire Category Theorem is not reversible, in that there are BaireCat spaces which are neither locally-compact Hausdorff, nor complete-metric; [FK, 1978].)

---

[‡] We use the term 'Baire set' in place of the more conventional, but ponderous, phrase "a set with the Property of Baire". Our usage differs from the way 'Baire set' is used in topological measure theory.

The Baire-field is likely to include most dynamical properties of interest to Ergodic theorists, as the following two models, due to Solovay 1970, show. There is a model of ZFC, set theory with Choice, in which all the "projective sets" (any finite number of projections and complements of Borel sets) are Baire-measurable. Indeed, if our Ergodic theorist is willing to forego Choice, there is a model of ZF in which *all* subsets of Polish spaces are Baire sets.

[‡] Density implies (and –when $\mathbb{A}$ is Polish– is equivalent to) topological transitivity of the action of $\Phi$ on $\mathbb{A}$. This latter is all that the proof requires; however, in applications one usually establishes a dense orbit.

PROOF OF ZERO-ONE LAW. Applying a homeomorphism $\varphi \in \Phi$ to the regular-picture yields

$$B \;=\; \varphi(B) \;=\; \varphi(U) \,\triangle\, \varphi(M) \,.$$

Since $\varphi$ is a homeomorphism, $\varphi(U)$ is regular-open and $\varphi(M)$ is meager. Then uniqueness of the regular-picture implies that both $U$ and $M$ are $\Phi$-invariant.

If $U$ is void, then $B = M$ and is consequently meager. Otherwise, $U$ intersects the hypothesized dense-orbit and so, by $\Phi$-invariance, $U$ is dense. Thus $U$ is residual and consequently so is $B$.  ◆

## Application to Group Actions

In classical Ergodic Theory, the role of "time" was usually played by $\mathbb{Z}$. Later, actions of more general topological groups $\Gamma$ (usually amenable) were studied. In this article, we will restrict ourselves to *countable discrete groups* $\Gamma$, although both the definitions and the results apply more generally. Although $\Gamma$ may be non-abelian, we write $\Gamma$ additively, with identity $0 \in \Gamma$. As names for "times", we use $\alpha, \beta, \gamma \in \Gamma$.

On our probability space $\mathbb{Y}$, an ***action*** of $\Gamma$ is a collection $\mathbf{T} = \{T^\gamma \mid \gamma \in \Gamma\}$ of $\mu$-mpts so that $T^{\beta+\alpha} \;=\; T^\beta \circ T^\alpha$ and the mapping $(\gamma, y) \mapsto T^\gamma(y)$ is a measurable map

$$(\Gamma \times Y, \text{Borel} \times \mathcal{Y}) \;\longrightarrow\; (Y, \mathcal{Y}) \,.$$

We call $T^\gamma$ "the transformation at time $\gamma$".

Each bi-mp map $\varphi: (Y, \mu) \to (Y', \mu')$ defines a conjugate action $\mathbf{R} \coloneqq \varphi \mathbf{T} \varphi^{-1}$ on $(Y', \mu')$ by $R^\gamma \coloneqq \varphi T^\gamma \varphi^{-1}$. Two actions $(\mathbf{T} : Y, \mu)$ and $(\mathbf{R} : Y', \mu')$ are ***isomorphic*** if they are conjugate. If the measure spaces are understood, we write $\mathbf{T} \cong \mathbf{R}$.

**The coarse or weak-operator topology.** On the set of measure-preserving transformations of $\mathbb{Y}$, the ***coarse topology*** is defined by net-convergence,

$$T_i \underset{i}{\rightarrowtail} R \quad \text{iff} \quad \mu\big(T_i^{-1}(E) \,\triangle\, R^{-1}(E)\big) \underset{i}{\rightarrowtail} 0, \text{ for all measurable subsets } E \subset Y.$$

A countable algebra $(E_k)_{k=1}^\infty$ of sets, which separates points of $Y$, gives rise to a metric which induces the coarse topology:

$$\text{dist}(T, R) \;\coloneqq\; \sum_{k=1}^\infty \frac{1}{2^k} \mu\big(T^{-1}(E_k) \,\triangle\, R^{-1}(E_k)\big) \,. \tag{1a}$$

With this metric, the set of mpts becomes a Polish space.

On $\mathbb{Y}$, let $\mathbb{A} = \mathbb{A}_\Gamma$ denote the set of $\Gamma$-actions. It also has a ***coarse topology*** from net-convergence: $\mathbf{T}_i \underset{i}{\rightarrowtail} \mathbf{R}$ iff $[\forall \gamma : T_i^\gamma \underset{i}{\rightarrowtail} R^\gamma]$. When $\Gamma$ is countable, the coarse topology makes the action space $\mathbb{A}$ a Polish space, via the metric

$$d(\mathbf{T}, \mathbf{R}) \;\coloneqq\; \sum_{j=1}^\infty \frac{1}{2^j} \text{dist}(T^{\gamma_j}, R^{\gamma_j}) \,, \tag{1b}$$

where $(\gamma_j)_{j=1}^\infty$ is a fixed enumeration of $\Gamma$.

**Weak Rohlin Property.** Zero-One applies to $\Gamma$ whenever $\Phi$ acts topologically transitively on $\mathbb{A}_\Gamma$. In the case $\Gamma = \mathbb{Z}$, the Rohlin Lemma implies this, since any two aperiodic transformations, $T$ and $R$, have arbitrarily large congruent Rohlin stacks. Thus some isomorphic copy $\varphi(T)$ is close to $R$, and so the set of copies, $\Phi(T)$, is dense among the aperiodics. The argument is finished by showing the aperiodics to be dense in $\mathbb{A}_\mathbb{Z}$.

With this as motivation, say that $\Gamma$ has the ***weak Rohlin Property*** (WRP) if $\Phi(\mathbf{T})$ is dense in $\mathbb{A}_\Gamma$, for some $\Gamma$-action $\mathbf{T}$. The analog of aperiodicity is "free-ness"; an action $\mathbf{T}$ is ***free*** if, after deleting a $Y$-nullset: *For each point $y$, the map $\gamma \mapsto T^\gamma(y)$ is injective.* The Appendix has a brief argument for the following:

$$\text{When } \Gamma \text{ is countable, the set of free } \Gamma\text{-actions is dense } -\text{indeed, residual}- \text{ in } \mathbb{A}_\Gamma. \tag{2}$$

Since there is a Rohlin lemma for $\Gamma = \mathbb{Z}^{\times d}$ (first shown in [KW] and later in [Con, thm 3.1]), the integer lattice has the WRP. More generally, the seminal paper by Ornstein & Weiss, [OW, p. 59], proves a multiple-tower Rohlin lemma. Weiss and Rudolph have indicated to us that this can be used to establish the WRP for all countable discrete amenable groups, thereby establishing that the Zero-One Law applies to these groups $\Gamma$.

*Question.* How large is the class of WRP groups; in particular, does it include the free group on two generators?

**Prime Example, 3.** A $\Gamma$-action $(\mathbf{T} : Y, \mathcal{Y}, \mu)$ is ***prime*** if it has no proper non-trivial factors. *Is primeness generic?* Even when $\Gamma = \mathbb{Z}$, this question is open, [dJR, P. 556]. What we do know, for weak-Rohlin groups $\Gamma$, is that Zero-One is applicable, since —as shown next— primeness is co-analytic.

Let $\mathcal{M}$ be the "measure-algebra" of $\mathcal{Y}$ (i.e, Borel sets equal a.e are identified) but without the null and co-null sets. Since $\mathcal{M}$ is Polish, $\mathcal{M} \times \mathcal{M}$ is Polish. Thus, to show that the set of <u>non</u>-prime actions, $\mathbf{N} \subset \mathbb{A}_\Gamma$, is analytic, we find a Borel condition (i.e, subset) $\mathbf{B}$ in $\mathcal{M} \times \mathcal{M} \times \mathbb{A}$ so that: *The projection of $\mathbf{B}$ on the $\mathbb{A}$-component, equals $\mathbf{N}$.*

And $\mathbf{N}$ *is* a projection, since it comprises those $\mathbf{T}$ for which: $\exists I, D \in \mathcal{M}$ such that *"$D$ is not in the $\mathbf{T},\mathcal{P}$-factor"*, where here —and henceforth— $\mathcal{P}$ abbreviates the two-set partition $(I, Y \smallsetminus I)$. Restating the quoted condition,

$$\exists \varepsilon : \ \forall \mathbf{T}, \mathcal{P}\text{-cylinder-sets } C : \ \mu(C \,\triangle\, D) \geq \varepsilon. \tag{*}$$

So $\mathbf{B}$ is the set of triples $(I, D, \mathbf{T})$ fulfilling $(*)$. It remains to show that $\mathbf{B}$ is Borel in $\mathcal{M} \times \mathcal{M} \times \mathbb{A}$.

*Cylinders.* Each finite subset $\mathcal{F} \subset \Gamma$ and tuple $\vec{v} = (v_1, \ldots, v_L)$ of maps $\mathcal{F} \to \{0, 1\}$ specifies a "$\mathbf{T}, \mathcal{P}$-cylinder set". It is $C_{\vec{v}}[I, \mathbf{T}] \coloneqq \bigcup_{\ell=1}^{L} C_{v_\ell}$, where $I_0 \coloneqq I$ and $I_1 \coloneqq Y \smallsetminus I$, and $C_v$ abbreviates $\bigcap_{\gamma \in \mathcal{F}} T^{-\gamma}(I_{v(\gamma)})$. Observe that the mass on the lefthand side of

$$\mu\big(C_{\vec{v}}[I, \mathbf{T}] \,\triangle\, D\big) \ \geq \ \varepsilon \tag{**}$$

varies continuously as a function of $I$ and $\mathbf{T}$ and $D$. So the set of $(**)$-triples $(I, D, \mathbf{T})$, is closed. Consequently, first intersecting over $\mathcal{F}$ and $\vec{v}$, then unioning over countably many $\varepsilon \searrow 0$, shows that $\mathbf{B}$ is an $\mathcal{F}_\sigma$-set. ◆

*Question.* Is primeness Borel measurable? We do not know. Allowing $\Gamma$ to range over continuous as well as discrete groups, the answer in principle could depend on $\Gamma$. For the related property of *simplicity*, [dJR], both its Borel and meager/residual status (the 0-1 law applies) are unknown.

## §2 Dynamically generically-equivalent settings

Two of the classical settings for studing $\Gamma$-actions were: *The action space,* $\mathbb{A}$, *where the measure was fixed and the transformation varied,* and *The set of shift-invariant measures, where the transformation was fixed, and the measure varied.* Our goal, the Equivalence Theorem[‡] below, is to show that these two settings give identical genericity results for Baire-measurable properties.

To discuss this generally, a **setting** is a topological space $\Omega$ and mapping $\omega \mapsto (\mathbf{R}_\omega : X_\omega, \mu_\omega)$ from $\Omega$ to $\Gamma$-actions. A subset $B \subset \Omega$ is **saturated** if $\mathbf{R}_\omega \cong \mathbf{R}_\eta$ implies that both $\omega$ and $\eta$, or neither, are in $B$. In principle, a saturated property might be Baire-measurable in one setting and not in another. Say that settings $\Omega$ and $\Omega'$ are **dynamically generically-equivalent** if saturated properties are Baire/meager/residual in $\Omega$ iff they are Baire/meager/residual in $\Omega'$.

Notes: Here is an example of a setting $\Omega$ where the Zero-One Law does not hold. In particular, it is a setting which is *not* dynamically generically-equivalent to $\mathbb{A}$.

Let $\Omega$ be the set of mpts on $(Y, \mathcal{Y}, \mu)$, but equipped with the "Halmos metric", rather than the (1a) metric. In the Halmos metric, $\mathrm{dist}(T, R)$ is the $Y$-mass of the set of points $y$ for which $T(y) \neq R(y)$. This metric makes $\Omega$ complete, but not separable.

Evidently, the open ball of radius $1/2$, centered at $Id$, is an isomorphism invariant open set, which is not residual. ♦

End of Notes.

### The space of Shift-Invariant Measures

We now specialize our measure space to be $\mathbb{Y} = (Y, \mathcal{Y}, \mathbf{m})$, where $Y := [0, 1)$ is topologized to be the unit circle, and $\mathbf{m}$ is Lebesgue measure (restricted to $\mathcal{Y}$, the Borel-field of the circle).

For our new setting, ${}^\bullet\mathbb{SIM}$, we need the definitions which follow. A superscript "$\bullet$" on a set of measures is meant to suggest that these measures might have atoms.

${}^\bullet\mathbb{M}$ is the set of (Borel probability) measures on $Y$.

$\mathbb{M}$ comprises those non-atomic $\mu \in {}^\bullet\mathbb{M}$ which have "full-support", that is, each non-void open set $V \subset Y$ has positive $\mu$-mass. We'll call these the **good** measures.

$\widehat{Y}$ denotes the product space $Y^{\times \Gamma}$, which is the infinite torus equipped with the product topology. Use $\widehat{y}, \widehat{z} \in \widehat{Y}$. Let $\widehat{y}[\gamma]$ denote the $\gamma^{\text{th}}$-coordinate of $\widehat{y}$.

$\mathbf{S}$ is the **shift** on $\widehat{Y}$. It is the $\Gamma$-action $\{S^\alpha \mid \alpha \in \Gamma\}$ defined by

$$S^\alpha(\widehat{y}) := \widehat{z}, \qquad \text{where} \quad \widehat{z}[\gamma] := \widehat{y}[\gamma + \alpha].$$

Note that $S^\beta \circ S^\alpha = S^{\beta + \alpha}$.

${}^\bullet\mathbb{SIM}$ is the set of shift-invariant measures on $\widehat{Y}$, that is, invariant under each $S^\alpha$, for $\alpha \in \Gamma$. We'll call such a measure a **sim**. For a sim $\widehat{\mu}$, let $\mathrm{Marg}[\widehat{\mu}]$ be the measure in ${}^\bullet\mathbb{M}$ which is the (1-dimensional) marginal of $\widehat{\mu}$; this marginal is independent of which copy of $Y$ we condition on, since $\widehat{\mu}$ is shift-invariant.

$\mathbb{SIM}$ is the set of **good sims**, that is, those $\widehat{\mu} \in {}^\bullet\mathbb{SIM}$ such that $\mathrm{Marg}[\widehat{\mu}]$ is in $\mathbb{M}$.

$\mathbb{SIM}_\nu$ comprises those good sims whose marginal is $\nu$, where $\nu$ is a good measure. Thus $\mathbb{SIM}$ is the disjoint union, over all $\nu \in \mathbb{M}$, of "fibers" $\mathbb{SIM}_\nu$.

---

[‡]In the case that $\Gamma = \mathbb{Z}$, essentially this same result was independently obtained by Dan Rudolph in his work [Rud] on orbit-equivalence.

**Weak-∗ preliminaries.** Both $^\bullet\mathbb{M}$ and $^\bullet\mathbb{SIM}$ become compact metric spaces under the weak-∗ topology: $\mu_i \to \nu$ iff $\left[\int t \, d\mu_i \longrightarrow \int t \, d\nu\right]$, for all continuous real-valued "test" functions $t$. For $U$ open, the indicator $\mathbf{1}_U$ can be arbitrarily well approximated by functions $t \leq \mathbf{1}_U$. Consequently, whenever $\mu_i \to \nu$,

$$\nu(\overline{B}) \;\geq\; \limsup_i \mu_i(\overline{B}) \;\geq\; \liminf_i \mu_i(B^\circ) \;\geq\; \nu(B^\circ) \tag{4}$$

for all sets $B$. In particular: *$\mu \mapsto \mu(B)$ is continuous on the collection of measures for which the boundary, $\partial B$, is a nullset.*

A *partition* $\mathcal{P}$ divides $Y$ (disjointly) into finitely many pieces; each *piece* $I \in \mathcal{P}$ is an interval of positive length. Fix a finite set $\mathcal{F} \subset \Gamma$ of times. An "*$\mathcal{F},\mathcal{P}$-cylinder set*", $C \subset \widehat{Y}$, is of the form $\bigotimes_{\gamma \in \mathcal{F}} I_\gamma$, where each $I_\gamma$ is in $\mathcal{P}$. This cylinder $C$ comprises those $\widehat{y}$ such that $\big[\forall \gamma \in \mathcal{F}\colon \widehat{y}[\gamma] \in I_\gamma\big]$. Observe that $\partial C$ is a nullset, for good sims $\widehat{\mu}$, since $\mathrm{Marg}[\widehat{\mu}]$ is non-atomic. Thus

On $\mathbb{SIM}$:     $\widehat{\mu}_i \to \widehat{\nu} \iff \forall\text{cylinders } C\colon \widehat{\mu}_i(C) \to \widehat{\nu}(C)$ $\tag{5}$

Implication ($\Leftarrow$) holds since, if sim $\widehat{\nu}$ differs from $\widehat{\eta}$, then $\widehat{\nu}(C) \neq \widehat{\eta}(C)$ for some $\mathcal{G}_\delta$-set $C$, hence for some open $C$, hence for some cylinder.

Finally, we leave to the Appendix a proof that

$$\mathbb{SIM} \text{ is a residual subset of } {}^\bullet\mathbb{SIM}. \tag{6}$$

Courtesy of this, *henceforth all measures and sims discussed are good*. Certainly $\mathbf{m}$ is good. We use names $\mu, \nu \in {}^\bullet\mathbb{M}$ for other good measures, and use $\widehat{\mu}, \widehat{\nu}, \widehat{\eta}$ to name good sims.

EQUIVALENCE THEOREM, 7. *Settings $\mathbb{A}$ and $^\bullet\mathbb{SIM}$ are dynamically generically-equivalent.*

**Overview of the proof.** Produce a BaireCat space $\mathbb{H}$ and an embedding $E$, which is a homeo­morphism from $\mathbb{H} \times \mathbb{A}$ into $\mathbb{SIM}$, so that

(8a) *For each $h \in \mathbb{H}$ and $\mathbf{T} \in \mathbb{A}$:  $E(h, \mathbf{T})$ is isomorphic to $\mathbf{T}$.*

(8b) *The set, Range($E$), is residual in $\mathbb{SIM}$.*

The map $E$ makes $\mathbb{H} \times \mathbb{A}$ into a setting. That $\mathbb{A}$ and $\mathbb{H} \times \mathbb{A}$ are dynamically generically-equivalent follows from (8a) and the plausible fact, [Oxt, P.57] theorems 15.2 and 15.3, that for arbitrary topological spaces with $\mathbb{H}$ a BaireCat space and $\mathbb{A}$ countably-generated,

*A subset $C \subset \mathbb{A}$ is Baire/meager/residual in $\mathbb{A}$  iff*
*the product $\mathbb{H} \times C$ is Baire/meager/residual in $\mathbb{H} \times \mathbb{A}$.*

Lastly, (8b) and (6) allow Proposition 9, below, to show that Range($E$) and $\mathbb{SIM}$ and $^\bullet\mathbb{SIM}$ are dynamically generically-equivalent.

**Relative Baire necessities.** Given subsets $B$ and $G$ of a topological space $\Omega$, say that "*$B$ is $G$-open*" (closed, nowhere-dense, meager, residual, Baire) iff $B \cap G$ is open (closed, etc.) relative to the induced topology on $G$.

PROPOSITION 9. *Suppose $M, B$ and $R$ are subsets of a topological space $\Omega$, with $R$ residual.*

(i) *If $M$ is $R$-meager then $M$ is meager. Also:  $B$ is $R$-Baire $\implies$ $B$ is Baire.*

(ii) *Suppose $\Omega$ is a BaireCat space. Then:*
    *$M$ is meager $\implies$ $M$ is $R$-meager;*
    *$B$ is Baire $\implies$ $B$ is $R$-Baire.*

Consequently, when $\Omega$ is BaireCat, a subset is Baire/meager/residual iff it is $R$-Baire/meager/residual.

PROOF OF (9$i$). Since the part of $M$ protruding from $R$ is meager, being inside of $\Omega \smallsetminus R$, we may assume that $M$ is included in $R$. Further, we may assume that $M$ is $R$-closed and has no $R$-interior.

Note that $\overline{R} \smallsetminus R$ has no interior. Our goal is to establish that $V \coloneqq \text{Interior}(\overline{M})$ is empty. But $V \subset \overline{M} \subset \overline{R}$. So it suffices to show that

$$V \text{ does not intersect } R,$$

for then $V \subset \overline{R} \smallsetminus R$ and is consequently void.

For the sake of contradiction, fix a point $\omega \in V \cap R$. Each neighborhood $W \ni \omega$ hits $M$ and so each $R$-neighborhood $W \cap R$ of $\omega$ also hits $M$, since $M \subset R$. Recall that $\omega \in R$ and $M$ is $R$-closed; thus $\omega \in M$.

This shows that $V \cap R \subset M$. But $M$ has no $R$-interior and so $V \cap R$ is empty, as desired.

Finally, for our $R$-Baire set $B$, we may assume that $B \subset R$, since $B \smallsetminus R$ is meager, hence Baire. So we can write $B = (V \cap R) \triangle M$, with $V$ open and with $M$ an $R$-meager set. From the foregoing, $M$ is meager, hence Baire. Since the Baire sets form an algebra, and $V, R, M$ are all Baire, $B$ must be a Baire set.                                                          ◆

PROOF OF (9$ii$). Our meager $M$ can be assumed to be nowhere-dense and closed; hence $M$ is $R$-closed. To show that $M$ is $R$–nowhere-dense, take an open $V$ so that $V \cap R$ is the $R$-interior of $M$. Now write $R = \Omega \smallsetminus M'$, with $M'$ meager. Thus $V \subset M \cup M'$. Since this union is meager, the openness of $V$ forces $V$ to be void. Consequently, $M \cap R$ is $R$-meager.

Finally, write Baire set $B$ as $U \triangle M$ with $U$ open and $M$ meager. Then

$$B \cap R \;=\; (U \cap R) \triangle (M \cap R) \;=\; R\text{-open} \triangle R\text{-meager}$$

and is consequently a Baire set relative to $R$.                              ◆

**Pushing measures forward.** Given a measurable map $f : (X, \mathcal{X}) \to (X', \mathcal{X}')$, each measure $\mu$ on $X$ pushes forward to a measure $f\langle \mu \rangle$ on $\mathcal{X}'$ defined by

$$f\langle \mu \rangle(B') \;\coloneqq\; \mu\big(f^{-1}(B')\big).$$

A (good) measure $\nu$ gives rise to a homeomorphism, $h$, of the circle $Y$, by

$$h(z) \coloneqq y, \text{ where } z = \nu\big([0, y)\big)$$

for $y \in Y$. Consequently, $\nu\big([0, y)\big) = \mathbf{m}\big([0, z)\big) = \mathbf{m}\big(h^{-1}([0, y))\big)$. Thus $h$ "adapts", so-to-speak, Lebesgue measure into $\nu$,

$$\nu = h\langle \mathbf{m} \rangle.$$

This $h$ is an order-preserving homeomorphism of $Y$ which fixes 0. We will call such a homeomorphism an **adaptation**, and let $\mathbb{H}$ denote the set of adaptations. Easily, $\mathbb{H}$ is a $\mathcal{G}_\delta$-subset of the Polish space $C(Y)$ of continuous functions—hence $\mathbb{H}$ is Polish [Oxt, pp. 47–50] and thus is BaireCat, as needed by **Overview**. Finally, note that

> *Each fiber* $\mathrm{SIM}_\nu$ *is the homeomorphic image of* $\mathrm{SIM}_{\mathbf{m}}$,
>
> *via the map* $\widehat{\mu} \mapsto \widehat{h}\langle\widehat{\mu}\rangle$;    (10a)

where here –and henceforth– for each adaptation $h$, we let $\widehat{h}$ denote the product $h^{\times\Gamma}$, which is a homeomorphism of $\widehat{Y}$.

*Pushing actions to sims.* Given an action $\mathbf{T} \in \mathbb{A}$, let $\widehat{\mathbf{T}}$ denote the map $y \mapsto \mathbf{T}(\cdot, y)$. This is an injection from $Y \hookrightarrow \widehat{Y}$ and gives a commutative diagram.

FIGURE 11   For each time $\gamma \in \Gamma$, the mapping $\widehat{\mathbf{T}}$ yields a commutative diagram. Lebesgue measure $\mathbf{m}$, being $T^\gamma$-invariant, pushes forward to an $S^\gamma$-invariant measure $\widehat{\mu} = \widehat{\mathbf{T}}\langle\mathbf{m}\rangle$. So $(\mathbf{S}\colon \widehat{Y}, \widehat{\mu})$, the action downstairs, is isomorphic to $(\mathbf{T}\colon Y, \mathbf{m})$.

$$
\begin{array}{ccc}
Y & \xrightarrow{\;T^\gamma\;} & Y \\
\Big\downarrow{\widehat{\mathbf{T}}} & & \Big\downarrow{\widehat{\mathbf{T}}} \\
\widehat{Y} & \xrightarrow{\;S^\gamma\;} & \widehat{Y}
\end{array}
$$

This mapping $\mathbf{T} \mapsto \widehat{\mathbf{T}}\langle\mathbf{m}\rangle$ is an embedding $\mathbb{A} \hookrightarrow \mathrm{SIM}_{\mathbf{m}}$. Its range is the set of **"graph-sims in $\mathrm{SIM}_{\mathbf{m}}$"**. As a first step towards residuality, we show that these are dense.

LEMMA 10b.  *The set of "graph-sims in* $\mathrm{SIM}_{\mathbf{m}}$ *" is dense in* $\mathrm{SIM}_{\mathbf{m}}$.

PROOF.  Fix a sim $\widehat{\nu}$ whose 1-marginal is $\mathbf{m}$. Given $\varepsilon$ and an $\mathcal{F},\mathcal{P}$-cylinder set, we wish to find an action $\mathbf{T} \in \mathbb{A}$ which gives $\varepsilon$-equality,

$$
\mathbf{m}\Big( \bigcap_{\gamma\in\mathcal{F}} \big[ T^{-\gamma} \big] (I_\gamma) \Big) \;\overset{\varepsilon}{\approx}\; \widehat{\nu}\Big( \bigotimes_{\gamma\in\mathcal{F}} I_\gamma \Big) .    (*)
$$

Fixing $\mathcal{P}$, it turns out that there is a $\mathbf{T}$ giving actual equality above, for all subsets $\mathcal{F} \subset \Gamma$ and all $\mathcal{F},\mathcal{P}$-cylinders.

For each piece $I \in \mathcal{P}$, let $\widehat{I}$ be corresponding time-0 set; the collection of $\widehat{y}$ for which $\widehat{y}[0] \in I$. Tautologically, $\widehat{\nu}(\widehat{I}) = \mathbf{m}(I)$ and the intersection $\bigcap_{\gamma\in\mathcal{F}}\big[S^{-\gamma}\big](\widehat{I}_\gamma)$ equals $\bigotimes_{\gamma\in\mathcal{F}} I_\gamma$.

Each $\widehat{I}$ is Borel-equivalent to an interval. Consequently there is a bi-mp Borel isomorphism $\varphi$ carrying probability space $(\widehat{Y}, \widehat{\nu})$ to $(Y, \mathbf{m})$, with $\varphi^{-1}(I) = \widehat{I}$ for each piece. Thus $\mathbf{T} := \varphi \mathbf{S} \varphi^{-1}$ gives equality in $(*)$.    ♦

**Measures which live on a graph.**  Temporarily, let $M$ denote the set of measures on $Y \times Y$ whose two marginals are equal and good. Consider an $\eta \in M$ satisfying, for all $B \in \mathcal{Y}$, that

$$
\exists A \in \mathcal{Y}\colon \quad \eta(A \times Y) \;=\; \eta(A \times B) \;=\; \eta(Y \times B) .    (12)
$$

Letting $\mu$ denote the marginal of $\eta$, this is equivalent to asserting that $\eta(A \times B)$ equals $\mu\big(A \cap R^{-1}(B)\big)$ for some $\mu$-mpt $R$. The graph of $R$ has full $\eta$-mass, and so $\eta$ is called a **"graph"** joining.

Let $\mathcal{C} \subset \mathcal{Y}$ be the countable algebra of sets formed by finite unions of rational-endpoint intervals in $Y$. It is straightforward to check that a measure $\eta$ is a graph joining iff (12) holds for each $B \in \mathcal{C}$. Equivalently, for $\varepsilon$ and all such $B$,

$$
\exists A \in \mathcal{C}\colon \quad \mathrm{Diameter}\big\{\eta(A \times Y), \eta(A \times B), \eta(Y \times B)\big\} \;<\; \varepsilon .
$$

Let $U_{\widehat{B}}^{\varepsilon}$ comprise such $\eta$. Since $\eta$ varies over measures with non-atomic marginals, the diameter above varies continuously with $\eta$. Thus $U_{\widehat{B}}^{\varepsilon}$ is open and, in consequence, the set $\bigcap_{B \in \mathcal{C}} \bigcap_{\varepsilon \searrow 0} U_{\widehat{B}}^{\varepsilon}$ –the graph joinings– is $\mathcal{G}_{\delta}$ in $M$.

*General graph-sims.* If each 2-dimensional marginal of sim $\widehat{\nu}$ is a graph joining, then $\widehat{\nu}$ is called a **graph-sim**. Letting $\nu$ be $\mathrm{Marg}[\widehat{\nu}]$, the shift-invariance of $\widehat{\nu}$ implies that it is of the form $\widehat{\nu} = \widehat{\mathbf{R}}\langle \nu \rangle$, for a (unique) action $(\mathbf{R} \colon Y, \nu)$.

Fix any two distinct times $\alpha, \beta \in \Gamma$ and consider the set of sims whose 2-marginal on $Y_{\alpha} \times Y_{\beta}$ is a graph. By the paragraphs above, this is a $\mathcal{G}_{\delta}$-set of sims. Since there are but countably many pairs in $\Gamma$, the set $\mathbb{SIM}^{\mathrm{Graph}}$ of graph-sims is $\mathcal{G}_{\delta}$ in $\mathbb{SIM}$. Consequently,

$$\mathbb{SIM}^{\mathrm{Graph}} \text{ is a residual subset of } \mathbb{SIM} \tag{10c}$$

since, courtesy of (10a,b), the graph-sims are dense in each fiber $\mathrm{SIM}_{\nu}$ and are thereby dense in $\mathbb{SIM}$.

---

## The embedding $E$

The preparation done, we can finally define embedding $E \colon \mathbb{H} \times \mathbb{A} \hookrightarrow \mathbb{SIM}$ which will fulfill the conditions of **Overview**. Let

$$E(h, \mathbf{T}) \;\coloneqq\; \widehat{h}\langle \widehat{\mathbf{T}}\langle \mathbf{m} \rangle \rangle \,, \tag{13}$$

where, as usual, $\widehat{h}$ means $h^{\times \Gamma}$. Evidently, the range of $E$ is simply $\mathbb{SIM}^{\mathrm{Graph}}$. In the next several paragraphs, we establish that $E$ is a homeomorphism by writing it as a composition of two homeomorphisms. The first one we consider is

$$D(h, \widehat{\mu}) \;\coloneqq\; \widehat{h}\langle \widehat{\mu} \rangle \,,$$

which is a mapping from $\mathbb{H} \times \mathrm{SIM}_{\mathbf{m}}$ onto $\mathbb{SIM}$.

**Is $D$ a homeomorphism?** Evidently $D$ is a bijection. In order to simplify the discussion of continuity properties of $D$, observe that for each adaptation $f \in \mathbb{H}$:

- *The composition map $h \mapsto fh$ is a homeomorphism of $\mathbb{H}$. Thus $\psi(h, \widehat{\mu}) \coloneqq (fh, \widehat{\mu})$ is a homeomorphism of $\mathbb{H} \times \mathrm{SIM}_{\mathbf{m}}$.*
- *$D$ carries $\psi$ to the mapping $\widehat{\eta} \mapsto \widehat{f}\langle \widehat{\eta} \rangle$, which is a homeomorphism of $\mathbb{SIM}$. (To see this, let $\widehat{\eta}$ be $\widehat{h}\langle \widehat{\mu} \rangle$.)*

Consequently, when arguing that $D$ is continuous at a point $(f, \widehat{\nu})$, or that $D^{-1}$ is continuous at $\widehat{f}\langle \widehat{\nu} \rangle$, we may assume that $f$ is the identity map $Id$ on $Y$. We now make these two arguments, having fixed: *A partition $\mathcal{P}$, a $\mathcal{P}$-cylinder set $C = \bigotimes_{\gamma \in \mathcal{F}} I_{\gamma}$, and a small $\varepsilon > 0$. Also, let $\delta$ denote $\varepsilon / |\mathcal{F}|$.*

*$D$ is continuous at $(Id, \widehat{\nu})$.* Fix $\widehat{\nu}$ and consider any $\widehat{\mu}$ sufficiently near $\widehat{\nu}$ that $\widehat{\mu}(C)$ is $\varepsilon$-close to $\widehat{\nu}(C)$. Take any adaptation $h$ for which

$$\left\| h - Id \right\|_{\mathrm{sup}} \;\leq\; \delta \,. \tag{$*$}$$

Then $\mathbf{m}\big(h^{-1}(I) \triangle I\big) \le 2\delta$, since each piece $I$ is an interval. Consequently,

$$\left| \widehat{\mu}\Big(\bigotimes_{\gamma \in \mathcal{F}} h^{-1}(I_\gamma)\Big) - \widehat{\mu}(C) \right| \le \sum_{\gamma \in \mathcal{F}} \mathbf{m}\big(h^{-1}(I_\gamma) \triangle I_\gamma\big)$$

$$\le |\mathcal{F}| \cdot 2\delta = 2\varepsilon. \qquad (**)$$

Applied to $C$, then, $\widehat{h}\langle\widehat{\mu}\rangle$ is $2\varepsilon$-close to $\widehat{\mu}$ and thus is $3\varepsilon$-close to $\widehat{\nu} = \widehat{Id}\langle\widehat{\nu}\rangle$. $\blacklozenge$

*The inverse, $D^{-1}$, is continuous at $\widehat{\nu} = \widehat{Id}\langle\widehat{\nu}\rangle$.* If $\widehat{h}\langle\widehat{\mu}\rangle$ is sufficiently near $\widehat{\nu}$, then their 1-marginals are close enough, choosing an $n$ with $\frac{1}{n} < \delta/2$, that

$$\mathbf{m}\Big(h^{-1}\big([0, y)\big)\Big) \overset{\delta/2}{\approx} \mathbf{m}\Big([0, y)\Big), \quad \text{for } y = \tfrac{1}{n}, \tfrac{2}{n}, \ldots, \tfrac{n}{n}.$$

Since $h$ preserves order on $Y$ and fixes zero, $h^{-1}(y)$ is $(\delta/2)$-close to $y$, for these points $y$. Thus $(*)$ holds. Therefore on $C$, as $(**)$ argued, sims $\widehat{h}\langle\widehat{\mu}\rangle$ and $\widehat{\mu}$ are $2\varepsilon$-close. So if $\widehat{h}\langle\widehat{\mu}\rangle$ is near enough to $\widehat{\nu}$ that they are $\varepsilon$-close on $C$, then $\widehat{\mu}$ is $3\varepsilon$-close (on $C$) to $\widehat{\nu}$.

Together with $(*)$, this establishes that the pair $(h, \widehat{\mu})$ is close to $(Id, \widehat{\nu})$. $\blacklozenge$

*Remark.* The two preceding arguments show that $D$ is a homeomorphism. By restricting $D$, we have the following:

*The map $(h, \widehat{\mu}) \mapsto \widehat{h}\langle\widehat{\mu}\rangle$ is a homeomorphism from $\mathbb{H} \times \mathrm{SIM}_{\mathbf{m}}^{\mathrm{Graph}}$ onto $\mathbb{SIM}^{\mathrm{Graph}}$,* (14)

where $\mathrm{SIM}_{\mathbf{m}}^{\mathrm{Graph}}$ is the set of graph-sims with marginal $\mathbf{m}$.

**A homeomorphism $\mathbb{A} \longrightarrow \mathrm{SIM}_{\mathbf{m}}^{\mathrm{Graph}}$.** Recall, from Figure 11, the injection $\mathbf{T} \mapsto \widehat{\mathbf{T}}\langle\mathbf{m}\rangle$. An argument similar to the continuity of $D$ shows that $\mathbf{T} \mapsto \widehat{\mathbf{T}}\langle\mathbf{m}\rangle$ is continuous. We now show that its inverse is continuous at an arbitrary point $\widehat{\mathbf{R}}\langle\mathbf{m}\rangle$, where $\mathbf{R} \in \mathbb{A}$.

Fix $\varepsilon$ and a time $\gamma$ and a subinterval $I \subset Y$. Let $J$ denote the subset $R^{-\gamma}(I)$. Regard $J \times I$ as a product set, on times $0$ and $\gamma$, and compute:

$$\widehat{\mathbf{R}}\langle\mathbf{m}\rangle(J \times I) = \mathbf{m}\big(J \cap R^{-\gamma}I\big) = \mathbf{m}(J);$$

$$\widehat{\mathbf{T}}\langle\mathbf{m}\rangle(J \times I) = \mathbf{m}\big(J \cap T^{-\gamma}I\big).$$

If sim $\widehat{\mathbf{T}}\langle\mathbf{m}\rangle$ is near enough to $\widehat{\mathbf{R}}\langle\mathbf{m}\rangle$, then the lefthand terms are $\varepsilon$-close.[†] So the righthand terms are $\varepsilon$-close. Hence $2\varepsilon > \mathbf{m}\big(T^{-\gamma}(I) \triangle J\big)$, which equals $\mathbf{m}\big(T^{-\gamma}(I) \triangle R^{-\gamma}(I)\big)$

This argument applies for all pieces $I$ in a partition, and for all times $\gamma$ in any specified finite set $\mathcal{F}$, so actions $\mathbf{T}$ and $\mathbf{R}$ have been shown to be close. $\blacklozenge$

**Final step.** The foregoing shows that the mapping $\big(h, \mathbf{T}\big) \mapsto \big(h, \widehat{\mathbf{T}}\langle\mathbf{m}\rangle\big)$ is a homeomorphism of $\mathbb{H} \times \mathbb{A}$ onto $\mathbb{H} \times \mathrm{SIM}_{\mathbf{m}}^{\mathrm{Graph}}$. Composing this with (14) shows that $E$ is a homeomorphism onto $\mathbb{SIM}^{\mathrm{Graph}}$, thus completing the proof of the *Equivalence Theorem.*

---

[†]This, by approximating $J$ with a finite union $J'$ of intervals, so that $\mathbf{m}(J' \triangle J) \le \varepsilon/2$.

## Acknowledgments

We thank our friends Karl Petersen –for the [Oxt2] reference– and Andrés del Junco, Dan Rudolph, and Benjy Weiss, for illuminating conversations on genericity issues. For $\mathbb{Z}$-actions, Dan put us wise to Lemma 10$b$, which he proved by a different argument. The second author gratefully thanks the University of Toronto for its hospitality while this paper was being completed.

## §A   Appendix

Proof of (6), first step: $\mathbb{M}$ *is a* $\mathcal{G}_\delta$-*subset of* $^\bullet\mathbb{M}$.   Write $\mathbb{M} = \mathbf{F} \cap \mathbf{N}$, where these denote the sets, respectively, of underline{full-support} and underline{non-atomic} measures in $^\bullet\mathbb{M}$. For two measures $\mu$ and $\nu$, the convolution is the measure

$$\big[\nu * \mu\big](A) \;:=\; \int_Y \mu(A \ominus y)\, d\nu(y)\,, \tag{15}$$

where $A \ominus y$ is the group translation $\{z \ominus y \mid z \in A\}$ in $Y$.

From (4), if $U \subset Y$ is open then $V_U^\varepsilon := \{\mu \in {}^\bullet\mathbb{M} \mid \mu(U) > \varepsilon\}$ is open. Letting $\mathcal{B}$ be a countable base for the topology of $Y$, then,

$$\mathbf{F} \;=\; \bigcap_{U \in \mathcal{B}} \bigcup_{\varepsilon \searrow 0} V_U^\varepsilon$$

is a $\mathcal{G}_\delta$-set.

For $C \subset Y$ closed, the set $\{\mu \in {}^\bullet\mathbb{M} \mid \mu(C) < \varepsilon\}$ is open, by (4). So $V_\varepsilon^K$ is open, where it comprises those $\mu$ such that, for each $n = 1, \ldots, K$:   $\mu\big(\big[\frac{n-1}{K}, \frac{n}{K}\big]\big) < \varepsilon$. Consequently,

$$\mathbf{N} \;=\; \bigcap_{\varepsilon \searrow 0} \bigcup_{K=1}^{\infty} V_\varepsilon^K$$

and so is $\mathcal{G}_\delta$ in $^\bullet\mathbb{M}$.

In order to show that $\mathbf{N}$ is dense in $^\bullet\mathbb{M}$, take $\delta$ a small positive number and let $\mathbf{m}_\delta := \frac{1}{\delta} \cdot \mathbf{m}|_{[0,\delta)}$ denote normalized Lebesgue measure on $[0, \delta)$. For each $\nu \in {}^\bullet\mathbb{M}$, then, the convolution $\nu * \mathbf{m}_\delta$ tends to $\nu$, as $\delta \searrow 0$. And $\nu * \mathbf{m}_\delta$ is non-atomic, since $\mathbf{m}_\delta$ is non-atomic. ◆

Second Step: $\mathbb{SIM}$ *is residual in* $^\bullet\mathbb{SIM}$.   Let $\mathbf{F}$ and $\mathbf{N}$ now denote the set of sims $\widehat{\mu} \in {}^\bullet\mathbb{SIM}$ whose time-0 marginal is a underline{full-support} or underline{non-atomic} measure. The foregoing shows that these are $\mathcal{G}_\delta$ sets in $^\bullet\mathbb{SIM}$. And since all 1-marginals of a sim are equal, $\mathbf{F} \cap \mathbf{N} = \mathbb{SIM}$.

To see that $\mathbf{F}$ is dense, take a particular sim $\widehat{\mu} \in \mathbf{F}$ (e.g, product measure $\mathbf{m}^{\times \Gamma}$). For each sim $\widehat{\nu}$, the average $(1 - \delta)\widehat{\nu} + \delta\widehat{\mu}$ is in $\mathbf{F}$, and tends to $\widehat{\nu}$ as $\delta \searrow 0$.

To show $\mathbf{N}$ dense in $^\bullet\mathbb{SIM}$, fix $\widehat{\nu}$ and $\delta > 0$. Again let $\mathbf{m}_\delta$ be normalized Lebesgue measure on $[0, \delta)$, and let $\widehat{\mathbf{m}_\delta}$ be the product measure $(\mathbf{m}_\delta)^{\times \Gamma}$. Taken over the torus $\widehat{Y}$, the convolution $\widehat{\nu} * \widehat{\mathbf{m}_\delta}$ is shift-invariant and has marginal $\mathrm{Marg}[\widehat{\nu}] * \mathbf{m}_\delta$, which is non-atomic. As $\delta \searrow 0$, furthermore, $\widehat{\mathbf{m}_\delta}$ tends to atomic measure on the identity $\widehat{0} \in \widehat{Y}$. Thus $\widehat{\nu} * \widehat{\mathbf{m}_\delta} \to \widehat{\nu}$. ◆

PROOF OF (2): *Free actions are dense in* $\mathbb{A}$. Let **Fr** be the set of "free sims" $\widehat{\mu}$, i.e, $(\mathbf{S} : \widehat{Y}, \widehat{\mu})$ is a free action. By showing **Fr** residual in $^\bullet\mathbb{SIM}$, the *Equivalence theorem* will tell us that the free actions are residual in $\mathbb{A}$.

Fix a time $\beta$ and let $C_\beta := \left\{ \widehat{y} \mid S^\beta(\widehat{y}) = \widehat{y} \right\}$. Since $S^\beta$ is continuous, $C_\beta$ is closed and consequently the collection $\left\{ \widehat{\mu} \in {}^\bullet\mathbb{SIM} \mid \widehat{\mu}(C_\beta) < \varepsilon \right\}$ is open. Intersecting over $\varepsilon$ and all times $\beta$ shows that **Fr** is $\mathcal{G}_\delta$. Certainly, **Fr** is dense in $^\bullet\mathbb{SIM}$ since, from above, $\widehat{\nu} * \widehat{\mathbf{m}_\delta} \to \widehat{\nu}$. And $\widehat{\nu} * \widehat{\mathbf{m}_\delta}$ is free. ◆